\documentclass[12pt]{article}
\usepackage{amssymb}
\usepackage{amsmath}
\numberwithin{equation}{section}

\newtheorem{theorem}{Theorem}

\newtheorem{corollary}[theorem]{Corollary}

\newtheorem{lemma}{Lemma}

\newcommand{\KKK}{K(\eta,y;x_0)}
\newcommand{\ooo}{\overline}
\newcommand{\ep}{\varepsilon}

\newcommand{\ppp}{\partial}
\newcommand{\www}{\widetilde}
\newcommand{\OOO}{\Omega}
\newcommand{\N}{\mathbb{N}}
\newcommand{\R}{\mathbb{R}}

\numberwithin{equation}{section}

\allowdisplaybreaks

\title
{\bf One-dimensional coefficient inverse problems by transformation operators}

\author{$^1$ Oleg Imanuvilov and  
$^{2,3}$ Masahiro Yamamoto}

\date{}
\begin{document}
\maketitle
\thanks{
$^1$ Department of Mathematics, Colorado State University, 101 Weber Building, 
Fort Collins CO 80523-1874, USA 
e-mail: {\tt oleg@math.colostate.edu}

$^2$ Graduate School of Mathematical Sciences, The University
of Tokyo, Komaba, Meguro, Tokyo 153-8914, Japan,\\
$^3$ Department of Mathematics, Faculty of Science, Zonguldak B\"ulent Ecevit 
University, \\
Zonguldak 67100, T\"urkiye\\
e-mail: {\tt myama@ms.u-tokyo.ac.jp}
}



\date{}
\maketitle

\begin{abstract}
We prove the uniqueness for an inverse problem of
determining a matrix  coefficient $P(x)$ of a system of  evolution equations
$\sigma \ppp_t u = \ppp_x^2 u(t,x) - P(x) u(t,x)$ for 
$0<x<\ell$ and $0<t<T$, where $\ell>0$ and $T>0$ are arbitrarily given.  
The uniqueness results assert that two solutions 
have the same Cauchy data at $x=0$ over $(0,T)$ and 
the same initial value or the final value which is positive on $[0,\ell]$,
then the zeroth-order coefficient is uniquely determined on $[0,\ell]$.
The uniqueness for inverse coefficient problem for a system of  
evolution equations without boundary conditions over the whole boundary 
is an open problem even in the one-dimension in the case where
only initial value is given as spatial data.
Moreover, in the case of the zero initial condition, we prove the 
uniqueness in the half of the spatial interval.
\end{abstract} 
\baselineskip 18pt

\section{Introduction and main results.}
Let $\ell, T>0$ and $N\in \N$.  Henceforth $\cdot^T$ denotes 
the transpose of vectors under consideration.

We consider a spatially one-dimensional linear evolution 
equation 
$$
L_{P}(x,D){u}:=\sigma\ppp_t{u} - \ppp_x^2{u} +P(x){u}=0,
\quad 0<x<\ell, \, 0<t<T,
$$ 
where ${u} = {u}(t,x)=(u_1(t,x),\dots, u_N(t,x))^T$ for 
$(t,x) \in (0,T) \times (0,\ell)$ and 
$P(x)=\{p_{ij}(x)\}_{1\le i,j \le N}$ is a coefficient depending 
on the spatial variable $x\in (0,\ell)$ and $\sigma$ is a nonzero 
complex number.
For $N=1$, we understand that $P(x)$ and $u(t,x)$ as scalar valued 
functions.

Typical examples of the operators $L_{P}(x,D)$ are as follows.
\\
{\bf (i) a single parabolic equation:}
$$
L_P(x,D)u:= \ppp_t u(t,x) - \ppp_x^2 u(t,x) + P(x) u(t,x), 
\quad 0<x<\ell,\, 0<t<T.
$$
{\bf (ii) a Schr\"odinger equation with $N=1$:}
$$
L_P(x,D) u:= \sqrt{-1}\ppp_tu(t,x) - \ppp_x^2u(t,x) + P(x)u(t,x), 
\quad 0<x<\ell,\, 0<t<T.
$$
{\bf (iii) a system of parabolic equations:}
$$
L_{P}(x,D){u}:= \ppp_t
\left( \begin{array}{c}
u_1(t,x) \\
\cdots \\
u_N(t,x) \\
\end{array}\right)
- \ppp_x^2
\left( \begin{array}{c}
u_1(t,x) \\
\cdots \\
u_N(t,x) \\
\end{array}\right)
+ P(x)\left( \begin{array}{c}
u_1(t,x) \\
\cdots \\
u_N(t,x) \\
\end{array}\right)
$$
for $0<x<\ell$ and $0<t<T$,
where $P(x) = (p_{ij}(x))_{1\le i,j \le N}$ is an 
$N\times N$ matrix.

Henceforth we consider $N$-number of $R^{N}$-valued functions
$$
u^k := \left( \begin{array}{c}
u_1^k(t,x) \\
\cdots \\
u_N^k(t,x) \\
\end{array}\right)
= (u_1^k, ..., u_N^k)^T \quad\mbox{for $k\in \{1,..., N\}$}.
$$
This article is first concerned with the following inverse coefficient 
problems:
\\
{\bf Inverse coefficient problem.}
{\it Let $T>0$ be arbitrarily fixed and for $k\in \{1, 2, ..., N\}$,
let $L_{P}(x,D){u}^k = 0$ for $(t,x) \in (0,T)\times (0, \ell) $
and $\ppp_x{u}^k (t,0) = 0$ for $0<t<T.$
Then determine $P(x)$, $0<x<\ell$ by data
$$
{ u}^k (0,\cdot)\vert_{(0,\ell)}, \quad { u}^k (\cdot,0)\vert_{(0,T)},
\quad 1\le k \le N
$$
or
$$
{u}^k (T,\cdot)\vert_{(0,\ell)}, \quad { u}^k (\cdot,0)\vert_{(0,T)},
\quad 1\le k \le N,
$$
provided that ${u}^k (0,\cdot)$ or ${u}^k (T,\cdot)$, $1\le k \le N$, are 
assumed to belong to some admissible set of functions. } 

In particular, for a single equation, 
our inverse problem can be formulated as follows:
\\
{\bf Inverse coefficient problem in the case of $N=1$.}
{\it Let $T>0$ be arbitrarily fixed and 
$L_P(x,D) = 0$ for $(t,x) \in (0,T)\times (0, \ell) $
and $\ppp_xu(t,0) = 0$ for $0<t<T$.
Then determine $P(x)$, $0<x<\ell$ by data
$$
u(0,\cdot)\vert_{(0,\ell)}, \quad u(\cdot,0)\vert_{(0,T)}
$$
or
$$
u(T,\cdot)\vert_{(0,\ell)}, \quad u(\cdot,0)\vert_{(0,T)},
$$
provided that $u(0,\cdot)$ or $u(T,\cdot)$ is assumed to belong to
some admissible set.}
\\

In this article, we propose a new methodology for proving the 
uniqueness for spatially one-dimensional inverse coefficient problems 
based on the transformation operator (e.g., Levitan \cite{Le}).
We state the uniqueness results for the three equations in the above
examples.

Henceforth we write $a'(x):= \frac{da}{dx}(x)$ for $0<x<\ell$, and 
$H^m(0,\ell)$ with $m\in \N$ denotes the usual Sobolev spaces.
Moreover we set 
$$
C^{1,2}( [0,T]\times  [0,\ell])
:= \{ {u}\in C([0,T]\times [0,\ell] );\, \ppp_t{u}, 
\ppp_x{u}, \ppp_x^2{u}
\in C([0,T]\times [0,\ell]) \}.
$$
We first state the uniqueness result for $N\in \N$:
\\
\begin{theorem} \label{Moon}
{\it
Let  $\sigma\in \Bbb C\setminus\{0\}$ and $P, Q \in (C^1[0,\ell])^{N\times N}$.
Let ${u}^k:= (u_1^k, ..., u_N^k)^T$, 
$\www{{u}}^k = (\www{u_1}^k, ..., \www{u_N}^k)^T
\in (C^{1,2}([0,T]\times [0,\ell]))^N$ with $1\le k \le N$ satisfy
\begin{equation}
\label{(1.6)}
\sigma\ppp_t{u}^k = \ppp_x^2{u}^k - P(x)
{u}^k\,\, \mbox{on}\,\,[0,T]\times [0,\ell], \quad 
\ppp_x{u}^k(t,0) = 0\, \, \mbox{on}\,\, (0,T),
\end{equation}
and
\begin{equation}\label{(1.7)}
\sigma\ppp_t\www{{u}}^k = \ppp_x^2\www{{u}}^k - {Q}(x)\www{{u}}^k\,\, 
\mbox{on}\,\,[0,T]\times [0,\ell] , \quad 
\ppp_x\www{{u}}^k(t,0) = 0\, \, \mbox{on}\,\, (0,T).
\end{equation}
\\
{\bf Case 1: spatial data at $t=0$.} 
\\
We assume 
\begin{equation}\label{(P1.8)}
\vert \mbox{det}\, ({u}^1(0,x) \cdots {u}^N(0,x))\vert 
> 0 \quad \mbox{on $[0, \ell]$}.                            
\end{equation}
Then equalities
\begin{equation}\label{(1.7)}
{u}^k(0,x) = \www{{u}}^k(0,x) \quad \mbox{on}\,\, (0,\ell)
\quad \mbox{for all}\,\,k\in  \{1,\dots, N\}
\end{equation} 
and
\begin{equation}\label{(1.seven)}
u^k(t,0) = \www{{u}}^k(t,0) \quad \mbox{for all $0<t<T$ and 
$k \in  \{1,\dots, N\}$}
\end{equation}
yield $P=Q$ on $[0,\ell]$.
\\
{\bf Case 2: spatial data at $t=T$.}
\\
We assume 
\begin{equation}\label{(1.88)}
\vert \mbox{det}\, ({u}^1(T,x) \cdots {u}^N(T,x))\vert 
> 0 \quad \mbox{on $[0, \ell]$}.                            
\end{equation}
Then (\ref{(1.seven)})
and 
\begin{equation}\label{(1.77)}
{u}^k(T,x) = \www{{u}}^k(T,x), \quad x\in (0,\ell)  
\quad \mbox{for all}\,\,k\in  \{1,\dots, N\}.                                   \end{equation}
imply $P={Q}$ on $[0,\ell]$.
}
\end{theorem}   

In Theorem \ref{Moon}, in order to determine $N^2$ coefficients, we are 
required to 
repeat observations by changing initial values ${u}^k(0,\cdot)$ or
final values $u^k(T,\cdot)$ with 
$1\le k \le N$ satisfying (\ref{(P1.8)}) or (1.6).
We remark that neither $P$ nor ${Q}$ is not necessarily 
symmetric.

Next we formulate our results for the case of single evolution 
equations.
\begin{equation}\label{(1.1)}
\mathcal{A} := \{ a\in H^3(0,\ell);\, a'(0) =  0, \quad
a>0 \,\,\mbox{on $[0,\ell]$}\}.                             
\end{equation}
\begin{corollary}\label{AX} (Single parabolic equation).
{\it
We assume that $N=1$, $\sigma = 1$, $P, Q \in C^1[0,\ell]$ are 
real-valued functions.
Let $u=u(t,x)$ and $\www{u} = \www{u}(t,x) 
\in C^{1,2}([0,T]\times [0,\ell])$
satisfy
\begin{equation}\label{(1.2)}
\left\{ \begin{array}{rl}
& \ppp_tu = \ppp_x^2u - P(x)u, \quad \ppp_xu(t,0) = 0,
\quad 0<x<\ell, \, 0<t<T,\\
& \ppp_t\www{u} = \ppp_x^2\www{u} - Q(x)\www{u}, \quad 
\ppp_x\www{u}(t,0) = 0, \quad 0<x<\ell, \, 0<t<T
\end{array}\right.
                                              \end{equation}
and
\begin{equation}\label{(1.3)}
u(0,\cdot) = \www{u}(0,\cdot) = a\in \mathcal{A} \quad \mbox{or} \quad
u(T,\cdot) = \www{u}(T,\cdot) = a\in \mathcal{A}.
                                                     \end{equation}
If 
$$
u(t,0) = \www{u}(t,0), \quad \mbox{on}\, \, (0,T),
$$
then $P=Q$ on $[0,\ell]$.
}
\end{corollary}

Corollary 2 is the uniqueness without any data at $x=\ell$ and 
is more desirable in practice,
and natural by recalling the unique continuation by Cauchy data at 
$x=0$ only (e.g., Mizohata \cite{Mi}, Saut and Schereur \cite{SS}), which can 
be stated for the one-dimensional case as follows.
If $u=u(t,x)$ satisfies 
$\ppp_tu = \ppp_x^2u - P(x)u$ for $0<x<\ell$ and $0<t<T$, and 
$u(t,0) = \ppp_xu(t,0) = 0$ for $0<t<T$, then 
$u=0$ in $(0,T)\times(0,\ell)$ without any boundary data at $x=\ell$.
Thus, comparing with the uniqueness in the Cauchy problem,
Corollary 2 provides a natural answer to this 
open uniqueness problem for the inverse problem in the one-dimensional case.

The uniqueness in Theorem \ref{Moon} and Corollary 2 is natural to be expected, but 
not known even for the one-dimensional case before.
Indeed, an effective method by Carleman estimates was established by 
Bukhgeim and Klibanov \cite{BK}, and the uniqueness was proved by data 
$$
u(t_0,\cdot)\vert_{(0,\ell)}, \,\, u(\cdot,0)\vert_{(0,T)}, \,\,
\ppp_xu(\cdot,0)\vert_{(0,T)}
$$
as long as we choose $t_0$ such that $0<t_0<T$.  As for such uniqueness in 
multidimensions, we can refer also to
Theorem 6.4.1 (p.152) in Isakov \cite{Is1}.
See also Imanuvilov and Yamamoto \cite{IY98},
Klibanov \cite{Kli}, Klibanov and Timonov \cite{KT},
Yamamoto \cite{Y09} as for more general and related results.
Here we do not intend any comprehensive list of references.
However, the method by Carleman estimate can not directly work for the 
case of $t_0=0$ or $t_0 = T$, so that the uniqueness in the cases  
$t_0=0$ and $t_0=T$ has been a severe open problem even for 
the one-dimensional parabolic equation, in spite of the significance.

In the case where the whole lateral boundary data are given, that is,
if $\ppp_xu(t,0)$ and $\ppp_xu(t,\ell)$ for $0<t<T$ are given,
we can refer to Imanuvilov and Yamamoto \cite{IY230}, Klibanov \cite{Kli}. 
The work \cite{Kli} considers parabolic equations in the whole space $\R^d$, 
and transforms the inverse parabolic problem into 
an inverse hyperbolic problem by means of an integral transform in 
time $t$ by assuming the time-analyticity of the solution to  
establish the uniqueness for the case $t_0=0$. 
Imanuvilov and Yamamoto \cite{IY230} proved 
the global Lipschitz stability for the case
of $t_0=T$, and the uniqueness for the case $t_0=0$ for 
parabolic equations in a bounded domain in $x$.
On the other hand, Imanuvilov and Yamamoto \cite{IY231} proved 
the uniqueness if an initial value
$a:= u(0,\cdot)$ belongs to a function space, which requires at least that 
$a \in \cap_{m=1}^{\infty} H^m(0,\ell)$.

Moreover we can refer to Murayama \cite{Mu}, Pierce \cite{Pi},
Suzuki \cite{Su}, Suzuki and Murayama \cite{SM}, where the 
uniqueness is established for a single parabolic equation with 
the full boundary conditions at $x=0$ and $x=\ell$.
The essence of their articles is a reduction of the inverse problem for a
parabolic equation to an inverse spectral problem, so that their uniqueness 
relies on the uniqueness in the inverse spectral problem and 
they applied the transformation operators.
For such a reduction, one must extract relevant spectral information from the 
solution to a parabolic equation, and so their works must assume 
that all the eigencomponents of an initial value are not zero, in other
words, their method does not work by assuming only the non-vanishing of 
initial value $a$ on 
$[0,\ell]$. We emphasize that our method does not need any results on the 
inverse spectral problems.
The method in Murayama \cite{Mu}, Pierce \cite{Pi}, 
Suzuki \cite{Su}, Suzuki and Murayama \cite{SM}
essentially depends on the reduction to 
the inverse eigenvalue problem using eigenfunction expansions of 
the solutions, and so can not work for the non-symmetric case, for 
example.  
Katchalov, Kurylev, Lassas and Mandache \cite{KKLM} proves 
the equivalence of the uniqueness for inverse problems for evolution 
equations and the corresponding inverse spectral problems, and 
comprehensively indicates an idea which was used in earlier
works \cite{Mu}, \cite{Pi}, \cite{Su}, \cite{SM}.
As for other works using the transformation operator, see 
Nakagiri \cite{Na1}, \cite{Na2}.

Furthermore, our method works for a wider class of one-dimensional 
evolution equations, and establishes the uniqueness in the case
$t_0=0$ and $t_0=T$, where the method by Carleman estimates 
can not work. 

From Theorem \ref{Moon}, we can directly derive

\begin{corollary}\label{AXX} 
{ (Schr\"odinger equation).}
{\it
Let $N=1$, $\sigma=\sqrt{-1}$, and $\mathcal A$ be defined by (\ref{(1.1)}).
We assume that $P, Q \in C^1[0,\ell]$.
Let $u=u(t,x)$ and $\www{u} = \www{u}(t,x) 
\in C^{1,2}([0,T]\times [0,\ell] )$
satisfy
\begin{equation}\label{(1.4)} 
\left\{ \begin{array}{rl}
& \sqrt{-1}\ppp_tu = \ppp_x^2u - P(x)u, \quad \ppp_xu(t,0) = 0,
\quad 0<x<\ell, \, 0<t<T,\\
& \sqrt{-1}\ppp_t\www{u} = \ppp_x^2\www{u} - Q(x)\www{u}, 
\quad \ppp_x\www{u}(t,0) = 0, \quad 0<x<\ell, \, 0<t<T.
\end{array}\right.
\end{equation}
Then 
$$
u(t,0) = \www{u}(t,0) \quad   \mbox{on}\, \, (0,T)
$$
and
$$
u(0,\cdot) = \www{u}(0,\cdot) =:a \in \mathcal{A}
$$
imply $P=Q$ on $[0,\ell]$.
}
\end{corollary}

Here we note that $P$, $Q$, $a$ are complex-valued.
The work by Baudouin and Puel \cite{BaP1}, \cite{BaP2} proves the 
global Lipschitz stability for multidimensions, under restrictive assumption 
\begin{equation}\label{(1.5)}
\mbox{$a$ or $\sqrt{-1}a$ is real-valued}        
\end{equation}
and the zero boundary conditions imposed on the whole boundary.
In Corollary \ref{AXX}, we do not need to assume  (\ref{(1.5)}).
As for related results on inverse problems for Schr\"odinger 
equations, see also Baudouin and Yamamoto \cite{BaY},
Bukhgeim \cite{Bu}, Imanuvilov and Yamamoto \cite{IYSch}, 
Yuan and Yamamoto \cite{YY}.

So far, we assume that $\vert u(0,\cdot)\vert \ne 0$ for 
$0\le x\le \ell$, but now
we will consider the case $u(0,\cdot) = 0$ in $(0, \ell)$
for single parabolic equations as follows:
\\
{\bf Inverse coefficient problem for a single parabolic equation
with zero initial condition.}
{\it Let $T>0$ be arbitrarily fixed and 
$L_P(x,D)u = 0$ for $(t,x) \in (0,T)\times (0, \ell) $
and $u(0,x) = 0$ on  $(0,\ell)$.
Then determine $P(x)$, $0<x<\ell$ by data
$$
\partial_xu(t,0)\vert_{(0,T)}, \quad u(t,0)\vert_{(0,T)}.
$$
}

We have 
\begin{theorem}\label{Moon1} 
{\it
We assume that $P, Q\in C^1[0,\ell]$ are real-valued, and 
$\sigma\in \Bbb C\setminus\{0\}$.
Let $u=u(t,x)$ and $\www{u} = \www{u}(t,x) 
\in C^{1,2}([0,T]\times [0,\ell] )$
satisfy
\begin{equation}\label{(1.44)} 
\left\{ \begin{array}{rl}
& \sigma\ppp_tu = \ppp_x^2u - P(x)u, \quad u(0,x) = 0,
\quad 0<x<\ell, \, 0<t<T,\\
& \sigma\ppp_t\www{u} = \ppp_x^2\www{u} - Q(x)\www{u}, 
\quad \www{u}(0,x) = 0, \quad 0<x<\ell, \, 0<t<T.
\end{array}\right.
\end{equation}
We assume that there exists some $m\in \Bbb N$ such that 
\begin{equation}\label{mobik}
\left\{ \begin{array}{rl}
& \ppp_t^mu, \ppp_t^m\www u \in C^{1,2}([0,T]\times [0,\ell]), \\
& \partial^m_t\partial_x u(0,0)\ne 0\quad\mbox{and} \quad 
\partial^j_t\partial_x u(0,0)=0\,\,\mbox{ for all}\,\,  0\le j \le m-1.
\end{array}\right.
\end{equation}
Then 
$$
u(t,0) = \www{u}(t,0),\quad \partial_xu(t,0) = \partial_x\www{u}(t,0)  \quad  
\mbox{on}\, \, (0,T)
$$
implies $P=Q$ on $\left[0,\frac{\ell}{2}\right]$.
}
\end{theorem}

We notice that the theorem asserts the uniqueness in the half interval 
of $[0,\ell]$ where we are given initial data.  The work Pierce \cite{Pi} is most 
related to Theorem \ref{Moon1} with the zero initial value, and we 
do not still require any boundary conditions at one end $x=\ell$.

{\bf Remark on zero initial values for the inverse problems}.
{\it The non-zero initial condition is important.  On the oher hand,
the articles \cite{Mu}, \cite{Su}, \cite{SM} assume 
not the positivity of an initial value $a$, but 
\begin{equation}\label{1.15}
\int^{\ell}_0 a(x)\psi_n(x) dx \ne 0
\end{equation}
for $n\in \N$, where 
$\{\psi_n\}_{n\in \N}$ is an orthonormal basis composed of 
the eigenfunctions of the operator $L_P(x,D)$ with the zero 
boundary condition of suitable type at $x=0, \ell$.  Although 
(\ref{1.15}) can be relaxed for $n\in \N$ except for a finite number of 
$n$, the condition (\ref{1.15}) is far away from a necessary condition.
Indeed, \cite{IY231} proves the uniqueness with a different  
assumption even in general dimensions, 
provided that $\vert a\vert > 0$ on $[0,\ell]$.  In particular, 
as a trivial corollary, the main result in 
\cite{IY231} implies the uniqueness in the case where (\ref{1.15}) 
holds only for a finite set of $n$, while for such a case, 
the works \cite{Mu}, \cite{Su}, \cite{SM} can not imply the uniqueness. 
On the other hand, we can conclude that non-zero initial values are 
essential for the uniqueness in the following sense.
More precisely, assumption (\ref{mobik}) is essential and can compensate for 
the zero initial value. 
Without this assumption, the uniqueness result of Theorem \ref{Moon1} fails.  
Indeed, let us fix smooth $P$ and $Q$ which are not identically equal 
on the segment $[0,\frac{\ell}{2}]$.
Let $u$ solve the initial boundary value problem 
$$
\ppp_tu = \ppp_x^2u - P(x)u, \quad u(0,x) = 0, \quad \partial_x u(t,0)=0, 
\quad  u(t,\ell)=g(t),
$$
where $g\in C^\infty[0,T]$, $g^{(j)}(0)=0$ for all $j\in\Bbb N$ and  
$g(t)>0$ on $(0,T)$. 
Consider $\widetilde u=u+ Ku$, where the operator 
$K$ is defined in Lemma \ref{main} below.
Then $\widetilde u$ solves the initial boundary value problem
$$
\ppp_t\www u = \ppp_x^2\www u - Q(x)\www u, \quad \www u(0,x) = 0, 
\quad \partial_x \www u(t,0)=0
$$
and 
$$
\www u(t,0)= u(t,0),
$$
but $P = Q$ does not necessarily follow.
}
\\

From Theorem \ref{Moon1}, we have the following result in determining 
two coefficients in the case $N=1$.

\begin{corollary}\label{Moon3} 
{\it
We assume that $r, \www r, P, Q \in C^1[0,\ell]$ and $\sigma\in \Bbb C
\setminus\{0\}$.
Moreover
\begin{equation}
r(0)=\www r(0).
\end{equation}
Let $u=u(t,x)$ and $\www{u} = \www{u}(t,x) 
\in C^{1,2}([0,T]\times [0,\ell] )$
satisfy
\begin{equation}\label{(1.44)} 
\left\{ \begin{array}{rl}
& \sigma\ppp_tu =  \ppp_x^2u-r(x)\partial_x u - P(x)u, \quad u(0,x) = 0,
\quad 0<x<\ell, \, 0<t<T,\\
& \sigma\ppp_t\www{u} = \ppp_x^2\www{u}-\www r(x)\partial_x \www{u} 
- Q(x)\www{u}, 
\quad \www{u}(0,x) = 0, \quad 0<x<\ell, \, 0<t<T.
\end{array}\right.
                                              \end{equation}
We assume that there exists $m\in \Bbb N$
\begin{equation}
\left\{ \begin{array}{rl}
& \ppp_t^mu, \, \ppp_t^m\www u \in C^{1,2}([0,T]\times [0,\ell]), \\
& \partial^m_t\partial_x u(0,0)\ne 0, \quad  
\partial^j_t\partial_x u(0,0)=0\,\,\mbox{ for all}\,\, 0\le j \le m-1.
\end{array}\right.
\end{equation}
Then 
\begin{equation}\label{sarkofag}
u(t,0) = \www{u}(t,0)\quad \mbox{and}\quad \partial_xu(t,0) 
= \partial_x\www{u}(t,0)  \quad   \mbox{on}\, \, (0,T)
\end{equation}
implies 
$$
P + \frac{1}{4}r^2 - \frac{1}{2}r'
= Q + \frac{1}{4}{\www r}^2 - \frac{1}{2}\www{r}' 
\quad \mbox{on $\left[0,\frac{\ell}{2}\right]$}.
$$
}
\end{corollary}

The assumption of the zero Neumann boundary condition at $x=0$
for unknown functions $\www u_k, u_k$ is restrictive.  We can  drop it 
instead by knowing the values of coefficients $P, Q$ near 
$x=0$. 
The result is given in the following theorem.
\begin{theorem} \label{AX} 
{\it
We assume that $P, {Q} \in (C^1[0,\ell])^{N\times N}$ and
$\sigma\in \Bbb C\setminus \{0\}$.
Let ${u}^k:= (u_1^k, ..., u_N^k)^T$, 
$\www{{u}}^k = (\www{u_1}^k, ..., \www{u_N}^k)^T
\in (C^{1,2}([0,T]\times [0,\ell]))^N$ with $1\le k \le N$ satisfy
\begin{equation}\label{(Xui1.2)}
\left\{ \begin{array}{rl}
&\sigma \ppp_t u^k = \ppp_x^2 u^k - P(x)u^k, 
\quad 0<x<\ell, \, 0<t<T,\\
&\sigma \ppp_t\www{u}^k = 
\ppp_x^2\www{u}^k - {Q}(x)\www{u}^k,  \quad 0<x<\ell, \, 0<t<T
\end{array}\right.
                                              \end{equation}
and either 
\begin{equation}\label{(P1.7)}
{u}^k(0,x) = \www{{u}}^k(0,x), \quad x\in (0,\ell) 
\quad \mbox{for all}\,\,k\in  \{1,\dots, N\}
\end{equation} or
\begin{equation}\label{(P1.77)}
{u}^k(T,x) = \www{{u}}^k(T,x), \quad x\in (0,\ell)  
\quad \mbox{for all}\,\,k\in  \{1,\dots, N\}.                                   \end{equation}
We assume 
\begin{equation}\label{(1.8)}
\vert \mbox{det}\, ({u}^1(0,x) \cdots {u}^N(0,x))\vert 
> 0 \quad \mbox{on $[0, \ell]$}                         
\end{equation}
in the case (1.21), and 
\begin{equation}\label{(1.88)}
\vert \mbox{det}\, ({u}^1(T,x) \cdots {u}^N(T,x))\vert 
> 0 \quad \mbox{on $[0, \ell]$}                            
\end{equation}
in the case (1.22).
Additionally suppose that there exists $\ep_0\in (0,\ell)$ such that 
$$ 
P(x)=Q(x)\quad \mbox{on}\quad (0,\ep_0).
$$                                                   
If 
$$
u^k(t,0) = \www{u}^k(t,0), \quad
\ppp_x u^k(t,0) = \ppp_x \www u^k(t,0) \quad \mbox{on}\, \, 
(0,T),
$$
then $P={Q}$ on $[0,\ell]$.
}
\end{theorem}

The article is composed of five sections.  
In Section \ref{two}, we show our main methodology and prove Theorem \ref{Moon}.
Sections \ref{three} and \ref{four} are
devoted to the proofs of Theorem \ref{Moon1} and Theorem \ref{AX} respectively.
In Section \ref{five}, we give concluding remarks.

\section{ Proof of Theorem \ref{Moon}.}\label{two}
Henceforth we set $\OOO := \{ (x,y);\, 0<y<x<\ell \}$.

First we show

\begin{lemma}\label{Asol}
{\it 
For $P, Q \in (C^1[0,\ell])^{N\times N}$, 
there exists a unique solution ${ K}=\{K_{ij}(x,y)\}_{1\le i,j \le N} 
\in (C^2(\ooo{\OOO}))^{N\times N}$ 
to the following problem:
\begin{equation}\label{(3.1)}
\left\{ \begin{array}{rl}
& \ppp_x^2{ K}(x,y) -\ppp_y^2 K(x,y)
= Q(x){ K}(x,y) - { K}(x,y)P(y), \quad (x,y) \in \OOO, \\
& \ppp_y {K}(x,0) = 0, \quad 0<x<\ell, \\
& 2\frac{d}{dx}{K}(x,x) = Q(x) - P(x),\quad 0<x<\ell, \quad
  K(0,0)=0
\end{array}\right.
\end{equation}
}
\end{lemma}

This is a Goursat problem and the proof is standard by means of the 
characteristics (e.g., Suzuki \cite{Su}).

We define an operator $K: (L^2(0,\ell))^N \longrightarrow 
(L^2(0,\ell))^N$ by 
$$
(Kv)(x) := \int^x_0  K(x,y)v(y) dy\quad 0<x<\ell.
$$
Here and henceforth, without fear of confusion, we use the same notation
for the operator $K$ with the integral kernel  $K(x,y)$.

We have

\begin{lemma}\label{main} 
Let $P, Q \in (C^1[0,\ell])^{N\times N}$, and 
let $u:= (u_1, ..., u_N)^T \in (C^{1,2}([0,T]\times [0,\ell]))^N$ 
satisfy 
\begin{equation}
\sigma\ppp_t{u} - \ppp_x^2{u} + P(x){u} = 0\quad\mbox{ in}\quad  (0,T)\times 
(0,\ell).
\end{equation} 
Then the function $\www{v}$ given by  
\begin{equation}\label{mk}
\www{v}(t,x) := u(t,x) + Ku(t,x) 
= u(t,x) + \int^x_0 K(x,y)u(t,y) dy, \quad 0<x<\ell, \, 0<t<T
\end{equation}
satisfies 
\begin{equation}\label{Mmk1}
\begin{split}
&\sigma\ppp_t\www{v} - \ppp_x^2\www{v} + Q(x)\www{v} 
= -K(x,0)\partial_x u(t,0)\quad\mbox{ in}\quad  (0,T)\times (0,\ell),\\ 
& \partial_x\www{v}(t,0) = \partial_x\www{u}(t,0),\quad 
\www{v}(t,0)={u}(t,0) \quad \mbox{for $0<t<T$}.
\end{split}\end{equation}
\end{lemma}

This is a classical transformation operator (e.g., 
Levitan \cite{Le}), and is used for one-dimensional inverse problems 
(e.g., Suzuki \cite{Su}, Suzuki and Murayama \cite{SM}).
We remark that unlike \cite{Su} and \cite{SM}, we do not assume the boundary 
value at $x=0$ for $u$, which produces the non-zero term on the right-hand
side of the first equation of (2.4).
Differently from \cite{Su}, \cite{SM} and \cite{Pi}, our method 
is free from any spectral properties of 
solutions, and so is applicable without full boundary conditions.

{\bf Proof.}
We directly verify that $\www{v}$ given by (\ref{mk}), satisfies 
(\ref{Mmk1}). 
Indeed,  
\begin{equation}\label{(3.2)}
\ppp_t\www{v}(t,x) = \ppp_t{u}(t,x) + \int^x_0 {K}(x,y)\ppp_t{u}(t,y) dy,
\quad 0<x<\ell, \, 0<t<T.                   
\end{equation}
Moreover, we have
$$
\ppp_x\www{v}(t,x) = \ppp_x{u}(t,x) +{K}(x,x){u}(t,x) 
+ \int^x_0 \ppp_x{K}(x,y){u}(t,y) dy.
$$
Therefore, (\ref{(3.1)}) yields 
\begin{equation}\label{(3.33)}
\begin{split}
& \ppp^2_x\int^x_0 {K}(x,y){u}(t,y) dy 
= \frac{d}{dx}( K(x,x) ){u}(t,x) + {K}(x,x)\ppp_x{u}(t,x)\\
+& (\ppp_x{K})(x,x){u}(t,x) + \int^x_0 \ppp_x^2K(x,y){u}(t,y) dy\\
= & \frac{d}{dx}({K}(x,x)){u}(t,x) + {K}(x,x)\ppp_x{u}(t,x)\\
+& (\ppp_x{K})(x,x){u}(t,x) + \biggl(\int^x_0 (\ppp_y^2 K)(x,y){u}(t,y) 
dy\\
+ & \int^x_0 (Q(x){K}(x,y)-{K}(x,y)P(y)){u}(t,y) dy\biggr), 
\quad 0<x<\ell, \, 0<t<T.
\end{split}
\end{equation}
Hence, in view of $\ppp_y{K}(x,0) =  0$ for $0<x<\ell$ which follows from 
the second equation in (\ref{(3.1)}), the integration by parts yields
\begin{equation}\label{(3.3)}
\begin{split}
&\int^x_0 (\ppp_y^2{K}(x,y)){u}(t,y) dy 
= \left[ (\ppp_y{K}(x,y)){u}(t,y)\right]^{y=x}_{y=0}
- \int^x_0 \ppp_yK(x,y)\ppp_y{u}(t,y) dy\\
=& (\ppp_yK)(x,x){u}(t,x) - \left[ {K}(x,y)\ppp_y{u}(t,y)\right]^{y=x}_{y=0}
+ \int^x_0 {K}(x,y)\ppp_y^2{u}(t,y) dy              \\
=& (\ppp_y{K})(x,x){u}(t,x) - {K}(x,x)\ppp_x{u}(t,x)\\ 
+ & K(x,0)\ppp_x{u}(t,0)
+ \int^x_0 {K}(x,y)\ppp_y^2{u}(t,y) dy\\
\end{split}
\end{equation}
By (\ref{(3.33)}) and (\ref{(3.3)}), we have
\begin{equation}\nonumber
\begin{split}
& \ppp^2_x\int^x_0 {K}(x,y){u}(t,y) dy 
  = \frac{d}{dx}({K}(x,x)){u}(t,x) + (\ppp_x{K})(x,x){u}(t,x) \\
+ & (\ppp_yK)(x,x)u(t,x) + {K}(x,0)\ppp_x{u}(t,0)
+ \int^x_0 {K}(x,y)\ppp_y^2{u}(t,y) dy\\
+ &\int^x_0 (Q(x){K}(x,y)-{K}(x,y)P(y)){u}(t,y) dy.
\end{split}
\end{equation}
Since
$$
\frac{d}{dx}({K}(x,x)) = (\ppp_x{K})(x,x) + (\ppp_y{K})(x,x),
$$ 
we can rewrite the above equality as
\begin{equation}\nonumber
\begin{split}
& \ppp^2_x\int^x_0 {K}(x,y){u}(t,y) dy 
= 2\frac{d}{dx}({K}(x,x)){u}(t,x)  \\
+ & \int^x_0 (Q(x){K}(x,y)-{K}(x,y)P(y)){u}(t,y) dy
+ {K}(x,0)\ppp_x{u}(t,0)\\
+ & \int^x_0 {K}(x,y)\ppp_y^2{u}(t,y) dy.
\end{split}
\end{equation}
Using the boundary condition (2.1), we obtain
\begin{equation}\label{saboromond}
 \ppp^2_x\int^x_0 {K}(x,y){u}(t,y) dy =
  Q(x){u}(t,x)-P(x){u}(t,x)  
\end{equation}
$$
+  \int^x_0 (Q(x){K}(x,y)-{K}(x,y)P(y)){u}(t,y) dy
+ {K}(x,0)\ppp_x{u}(t,0)
$$
$$
+ \int^x_0 {K}(x,y)\ppp_y^2{u}(t,y) dy, \quad 0<t<T,\, 0<x<\ell.
$$

Therefore, using $\sigma\ppp_t{u} - \ppp_x^2{u} + P(x){u} = 0$ 
in $(0,T)\times (0,\ell)$, we obtain
\begin{align*}
& \sigma \ppp_t\www{v} - \ppp_x^2\www{v} + Q(x)\www{v}
=\sigma \ppp_t{u} + \int^x_0 {K}(x,y)\sigma\ppp_t{u}(t,y) dy
- \ppp_x^2{u} \\
-& Q(x){u}+P(x){u}
- \int^x_0 (Q(x){K}(x,y)-{K}(x,y)P(y)){u}(t,y) dy 
- \int^x_0 K(x,y)\ppp_y^2{u}(t,y) dy            \\
+ & Q(x){{u}}+Q(x)\int^x_0 {K}(x,y){u}(t,y) dy - {K}(x,0)\ppp_x{u}(t,0)
= -K(x,0)\ppp_x{u}(t,0).
\end{align*}
Thus the proof of the lemma is complete.
$\blacksquare$
\\

Let $K(x,y) $ be determined by (\ref{(3.1)}).  
We now show a key lemma. 
\\
{\bf Lemma 3}
{\it
Let $u^k$ and $\www{u}^k$ satisfy (\ref{(1.6)}) and (\ref{(1.7)}).
For $x_0 \in (0, \ell)$ and $\delta\in (0, \ell-x_0)$, we assume
\begin{equation}\label{2.9}
K(x,0)\ppp_xu^k(t,0) = 0, \quad 0<x<x_0+\delta, \, 0<t<T, 
\, k\in \{1, ..., N\} 
\end{equation}
and
$$
P(x) = Q(x), \quad 0<x<x_0.
$$
Moreover we assume either (\ref{(P1.8)})-(\ref{(1.7)})-(\ref{(1.seven)}) or (\ref{(1.seven)})-(1.6)-(1.7).               Then there exists a constant $\ep > 0$ such that
$P=Q$ in $(x_0, \, x_0+\ep)$.
}
\\

Once Lemma 3 is proved, we can readily complete the proof of Theorem \ref{Moon} 
as follows.  It suffices to prove in the case of (\ref{(P1.8)})-(\ref{(1.7)})-(\ref{(1.seven)}), because the
case (\ref{(1.seven)})-(1.6)-(\ref{(1.77)}) is similar.
By $u^k, \www{u}^k \in C^{1,2}([0,T]\times [0,\ell])$, substituting 
$t=0$ and $x=0$ in (\ref{(1.6)}) and (1.2), and using $\www{u}^k(0,0) > 0$, we see
$P(0) = Q(0)$.  Therefore, we can define $x_0 \in [0,\ell]$ by the 
maximal point such that 
$$
P(x)=Q(x) \quad \mbox{for $0\le x \le x_0$}.
$$
If $x_0 = \ell$, then the proof is already finished, and so we can 
assume $0\le x_0 < \ell$.
By the assumption $\ppp_xu^k(t,0) = 0$ for $k\in \{1, ..., N\}$ and 
$0<t<T$, condition (2.9) is satisfied for $0<x<\ell$.  
Therefore Lemma 3 implies the 
existence of $\ep>0$ such that $P=Q$ in $(0, x_0+\ep)$.  This contradicts
the maximality of $x_0$.  Thus the proof of Theorem \ref{Moon} is complete.
$\blacksquare$

Now we proceed to 
\\
{\bf Proof of Lemma 3.}
Since $K(x,x) = 0$ for $0\le x\le x_0$ by the choice of $x_0$, 
the uniqueness of solution to (\ref{(3.1)}) yields 
\begin{equation}\label{XUK}
K(x,y)=0\quad \mbox{if $0 \le y \le x \le x_0$}.
\end{equation}

We define the functions $ \www{v}^k$ by
$$
\www{v}^k(t,x) = {{u}}^k(t,x) + \int^x_0 K(x,y)u^k(t,y)dy, \quad 
0<x<\ell, \, 0<t<T \quad \mbox{for all $k\in \{1,\dots,N\}$}.
$$
Then, in terms of $\ppp_xu^k(t,0) = 0$ for $0<t<T$, Lemma 2 implies
\begin{equation}\label{mk1}
\left\{ \begin{array}{rl}
& \sigma\ppp_t\www{v}^k - \ppp_x^2\www{v}^k + Q\www{v}^k = 0\,
\mbox{ in}\, (0,T)\times (0,\ell),\\
& \partial_x\www{v}^k(t,0)=0,\, \www{v}^k(t,0)={u}^k(t,0) \quad \mbox{on}\,\, 
[0,T].
\end{array}\right.
\end{equation}
By (\ref{(1.7)}), (\ref{(1.seven)}) and (\ref{mk1}), the function $w^k:= \www {u}^k- \www{v}^k$ 
satisfies 
$$
\sigma\partial_t w^k - \partial^2_x w^k+Q w^k=0 \,\mbox{ in}\,\, 
(0,T)\times (0,\ell),\quad w^k(t,0)=\partial_{x}w^k(t,0)=0,
\quad 0\le t \le T.
$$
By the classical uniqueness result (e.g., \cite{Mi}, \cite{SS}) 
of the Cauchy problem for the parabolic equation, we have 
$$
\www{v}^k=\www u^k \quad \mbox{in $(0,T)\times (0,\, \ell)$ 
for $k\in\{1,\dots,N\}$}.
$$
These equalities and the definition of $\www v^k$ by the transformation 
operator imply 
\begin{equation}\label{main1}
\int_0^xK(x,y) u^k(0,y)dy=0\quad \mbox{on}\quad  [0,\ell]\quad 
\mbox{for all $k\in \{1,\dots,N\}$}.
\end{equation} 
if (\ref{(1.7)}) is assumed, and 
\begin{equation}\label{Amain1}
\int_0^xK(x,y) u^k(T,y)dy=0\quad \mbox{on}\quad  [0,\ell]\quad 
\mbox{for all $k\in \{1,\dots,N\}$},
\end{equation}
if (\ref{(1.77)}) holds true.

Let (\ref{main1}) hold. The proof for the case (\ref{Amain1}) 
is the same. 

Using (\ref{saboromond}) and twice differentiating (\ref{main1})
with respect to $x$, we have
\begin{equation}\label{saboromond1}
\begin{split}
&0= \ppp^2_x\int^x_0 {K}(x,y) u^k(0,y)dy =
  Q(x) u^k(0,x)-P(x) u^k(0,x)  \\
+ & \int^x_0 (Q(x){K}(x,y)-{K}(x,y)P(y)) u^k(0,y) dy\\
+ & \int^x_0 {K}(x,y)\ppp_y^2 u^k(0,y) dy,
\quad 0<x<\ell.
\end{split}
\end{equation}

This equality implies
\begin{equation}\label{BBp}
\begin{split}
&0= 
  (Q(x)-P(x))B(x)  
+  \int^x_0 (Q(x){K}(x,y)-{K}(x,y)P(y))B(y) dy\\
+& \int^x_0 {K}(x,y)\ppp_y^2 B(y) dy, \quad 0<x<\ell
\end{split}
\end{equation}
where 
$$
B(x):= (u^1(0,x),\dots,u^N(0,x))
$$
is an $N\times N$ matrix. 
Assumption (\ref{(P1.8)}) implies 
$$
\vert \mbox{det} B(x)\vert>0, \quad 0\le x\le \ell.
$$
Applying the matrix $B^{-1}$ to both sides of equation (\ref{BBp}), 
we obtain
\begin{equation}\label{BBq}
\begin{split}
&0= 
  (Q(x)-P(x)) 
+   \int^x_0 (Q(x){K}(x,y)-{K}(x,y)P(y))B(y) dy\,B^{-1}(x)\\
+ &\int^x_0 {K}(x,y)\ppp_y^2 B(y) dy\, B^{-1}(x),
\quad 0<x<\ell.
\end{split}
\end{equation}
On the other hand, dividing the integral interval into 
$(0,x_0)$ and $(x_0, x)$, we write (\ref{BBq}) in the form
\begin{equation}\label{BBqq}
\begin{split}
&0= (Q(x)-P(x)) 
+   \int^x_{x_0} (Q(x){K}(x,y)-{K}(x,y)P(y))B(y) dy\,B^{-1}(x)\\
+ &\int^x_{x_0} {K}(x,y)\ppp_y^2 B(y) dy\, B^{-1}(x) \\
+ & \int^{x_0}_0 (Q(x){K}(x,y)-{K}(x,y)P(y))B(y) dy\,B^{-1}(x)\\
+ & \int^{x_0}_0 {K}(x,y)\ppp_y^2 B(y) dy\, B^{-1}(x), \quad
0<x<\ell
\end{split}
\end{equation}
Henceforth we set $\OOO_x= \{ (\xi, \eta);\, 0\le \eta \le \xi \le x\}$
for $x \in (0,\ell)$.

Estimating the second term on the right hand side of (\ref{BBqq}), we have
\begin{equation}\label{MANa}
\begin{split}
&\left\vert \int^x_{x_0} (Q(x){K}(x,y)-{K}(x,y)P(y))B(y) dy\,B^{-1}(x)
+ \int^x_{x_0} {K}(x,y)\ppp_y^2 B(y) dy\, B^ {-1}(x)\right\vert \\
\le & \vert x-x_0\vert\Vert B^{-1}\Vert_{C[0,\ell]} \Vert B\Vert_{C[0,\ell]}
 (\Vert P\Vert_{C[0,\ell]}+\Vert Q \Vert_{C[0,\ell]})
\sup_{(\xi, \eta) \in \ooo{\OOO_x}} \vert K(\xi, \eta)\vert   \\
+ & \vert x-x_0\vert\Vert B^{-1}\Vert_{C[0,\ell]} \Vert B\Vert_{C^2[0,\ell]}
\sup_{(\xi,\eta) \in \ooo{\OOO_x}} \vert K(\xi, \eta)\vert.
\end{split}\end{equation}

By the estimate of the solution $K$ to the Goursat problem (\ref{(3.1)}) 
for any $x\in [0,\ell]$,
we have
\begin{equation}\label{lox3}
\Vert K\Vert_{C^1(\ooo{\OOO_x})} \le C \Vert P-Q\Vert_{C[0,x]}.
\end{equation}
From (\ref{MANa}) and (2.19), we obtain
\begin{equation}\label{MANa1}
\begin{split}
&\left\vert \int^x_{x_0} (Q(x){K}(x,y)-{K}(x,y)P(y))B(y) dy\,B^{-1}(x)
+ \int^x_{x_0} {K}(x,y)\ppp_y^2 B(y) dy\, B^{-1}(x)\right\vert \\
\le & C\vert x-x_0\vert \Vert B^{-1}\Vert_{C[0,\ell]} \Vert B\Vert_{C[0,\ell]}
(\Vert P\Vert_{C[0,\ell]}+\Vert Q\Vert_{C[0,\ell]}) \Vert P-Q\Vert_{C[0,x]}\\
+& C\vert x-x_0\vert\Vert B^{-1}\Vert_{C[0,\ell]} \Vert B\Vert_{C^2[0,\ell]}
\Vert P-Q\Vert_{C[0,x]}\\
\le & C\vert x-x_0\vert \Vert P-Q\Vert_{C[0,x]}.
\end{split}
\end{equation}
Here and henceforth $C>0$ denotes generic constants depending on 
$B, P, Q, \ell$.
 
Now we will estimate the fourth and the fifth terms on the right-hand side 
of (\ref{BBqq}). Using (\ref{XUK}), 
we can write this term as
\begin{equation}\begin{split}
&I:= \int^{x_0}_0 (Q(x){K}(x,y)-{K}(x,y)P(y)) B(y) dy\,B^{-1}(x)
+ \int^{x_0}_0 {K}(x,y)\ppp_y^2 B(y) dy\, B^{-1}(x)\\
=& \int^{x_0}_0 \{ Q(x)({K}(x,y)-K(x_0,y))-({K}(x,y)-K(x_0,y))P(y)\}
 B(y) dy\,B^{-1}(x)\\
+ & \int^{x_0}_0 ({K}(x,y)-K(x_0,y))\ppp_y^2 B(y) dy\, B^{-1}(x),
\quad x\in [0, \ell].
\end{split}
\end{equation}
Then, applying mean value theorem, we obtain
\begin{equation}\label{xruk}\begin{split}
\vert I\vert\le x_0 \Vert B^{-1}\Vert_{C[0,\ell]} 
\Vert B\Vert_{C[0,\ell]} (\Vert P\Vert_{C[0,\ell]}+\Vert Q\Vert_{C[0,\ell]})
\sup_{y\in [0,x_0]}\vert {K}(x,y)-K(x_0,y)\vert   \\
+ x_0 \Vert B^{-1}\Vert_{C[0,\ell]} \Vert B\Vert_{C^2[0,\ell]}
\sup_{y\in [0,x_0]}\vert {K}(x,y)-K(x_0,y)\vert\\
\le C\sup_{(z,y)\in [0,x]\times [0,x_0]}\vert \partial_z K(z,y)\vert 
\vert x-x_0\vert.
\end{split}
\end{equation}
From (\ref{xruk}), using (2.19), we obtain
\begin{equation}\label{lox21}
\vert I\vert
\le C\Vert P-Q\Vert_{C[0,x]} \vert x-x_0\vert.
\end{equation}
From (2.17), (2.20) and (2.23), we have
\begin{equation}\nonumber
\vert (P-Q)(x)\vert \le C\vert x-x_0\vert \Vert P-Q\Vert_{C[0,x]}.
\end{equation}
This inequality implies
\begin{equation}\nonumber
\Vert P-Q\Vert_{C[0,x]}\le C \vert x-x_0\vert \Vert P-Q\Vert_{C[0,x]}
\end{equation}
for $0 < x < \ell$.
Taking sufficiently small $\ep>0$ satisfying $C\ep <1$, we obtain 
$$
\Vert P-Q\Vert_{C[0,x_0+\ep]}=0.
$$
This contradicts the choice of  $x_0$.
Thus the proof of Lemma 3 is complete. 
$\blacksquare$

\section{Proof of Theorem \ref{Moon1}}\label{three}

Let $u, v$ satisfy the conditions in Theorem \ref{Moon1} and
the function $\www v$ be determined by the formula (\ref{mk}):
$$
\www v(t,x) = u(t,x) + \int^x_0 K(x,y)u(t,y) dy.
$$
By Lemma \ref{main}, the function $\www v$ satisfies
\begin{equation}\label{PMmk1}\begin{split}
&\sigma\ppp_t\www{v} - \ppp_x^2\www{v} + Q(x)\www{v} 
= -\partial_x u(t,0) K(x,0)\quad\mbox{ in}\quad  (0,T)\times (0,\ell),\\ 
& \partial_x\www{v}(t,0)=\partial_x{u}(t,0), \quad \www{v}(t,0)={u}(t,0)\,\,
\mbox{on} \,\,(0,T),\quad \www{v}(0,\cdot)=0 \quad \mbox{in $(0,\ell)$}.
\end{split}
\end{equation} 
Moreover by $\ppp_t^mu\in C^{1,2}([0,T]\times [0,\ell])$,
we see that 
$\ppp_t^m \www v\in C^{1,2}([0,T]\times [0,\ell])$.
 
Then, the function $w=\www v-\www u$ satisfies 
$\ppp_t^m w\in C^{1,2}([0,T]\times [0,\ell])$ and  
\begin{equation}\label{UX}
\begin{split}
&\sigma\ppp_t w - \ppp_x^2 w + Q(x) w = -\partial_x u(t,0) K(x,0)
\quad\mbox{ in}\quad  (0,T)\times (0,\ell),\\ 
&\partial_x w(t,0) = w(t,0)=0\quad \mbox{in $(0,T)$}, \quad
w(0,\cdot)=0 \quad \mbox{in $(0,\ell)$}.
\end{split}\end{equation} 
Differentiating  equation (\ref{UX}) $m$-times with respect to $t$,
setting $\www w:= \ppp_t^mw \in C^{1,2}([0,T]\times [0,\ell])$, we obtain
\begin{equation}\label{UXX}
\begin{split}
&\sigma\ppp_t \www w - \ppp_x^2\www w + Q(x) \www w 
= -\partial_t^m \partial_xu(t,0) K(x,0)\quad\mbox{ in}\quad 
(0,T)\times (0,\ell),\\ 
& \partial_x\www w(t,0) = \www w(t,0) =0 \quad \mbox{in $(0,T)$},
\quad\www w(0,\cdot)=0 \quad \mbox{in $(0,\ell)$}.
\end{split}\end{equation} 
Here we can verify $\www w(0,\cdot) = 0$ as follows.
By (3.2) and (1.14), using $w(0,x) = 0$ for $0<x<\ell$, we see
$$
\sigma \ppp_tw(0,x) = - \ppp_xu(0,0)K(x,0) = 0, 
$$
that is, $\ppp_tw(0,x) = 0$ for $0<x<\ell$.  Differentiating 
(3.2) with respect to $t$ and substituting $t=0$, we obtain
$$
\sigma \ppp_t^2 w(0,x) = \ppp_x^2\ppp_tw(0,x) - Q(x)\ppp_tw(0,x)
- \ppp_x\ppp_tu(t,0)K(x,0).
$$
Therefore, $\ppp_t^2w(0,x) = 0$ for $0<x<\ell$
by (1.14) and $\ppp_tw(0,x) = 0$ for $0<x<\ell$.  In view of (1.14),
continuing this argument, we reach $\ppp_t^mw(0,x) = 0$, that is,
$\www w(0,x) =  0$ for $0<x<\ell$.
\\

We set
$$
R(t):= \ppp_t^m\ppp_xu(t,0), \quad 0<t<T,
$$
and
$$
(Mv)(t):= R(0) v(t) + \int^t_0 R'(t-\tau)v(\tau) d\tau, \quad 0<t<T.
$$
We recall $R'(t) = \frac{dR}{dt}(t)$.

We consider the equation $(Mz)(t,x) = \ppp_t\www w(t,x)$ for 
$0<t<T$ and $0<x<\ell$, that is,
\begin{equation}\label{3.4}
\ppp_t\www w(t,x) = R(0)z(t,x)
+ \int^t_0 R'(t-\tau)z(\tau,x) d\tau,
\quad 0<t<T,\, 0<x<\ell.
\end{equation}

Since $R(0) = \ppp_t^m\ppp_xu(0,0) \ne 0$ by (\ref{mobik}), 
the operator $M$ is a Volterra operator of the second 
kind.  Moreover, since $\int^t_0 R'(t-\tau)v(\tau)d\tau 
= \int^t_0 R'(\tau)v(t-\tau) d\tau$, we see that $M^{-1}: H^1(0,T) 
\longrightarrow H^1(0,T)$ exists and 
is bounded.  Therefore, $z \in H^{1,2}((0,T)\times (0,\ell))$ is well 
defined for each $x \in (0,\ell)$ by means of $\ppp_t\www w \in 
H^{1,2}((0,T)\times (0,\ell))$.

Since 
$$
R(0)z(t,x) + \int^t_0 R'(t-\tau)z(\tau,x) ds
= \ppp_t\left( \int^t_0 R(t-\tau)z(\tau,x) d\tau\right), 
$$
in view of (3.4), we have
$$
\ppp_t\left( \www w (t,x) - \int^t_0 R(t-\tau)z(\tau,x) d\tau\right)
= 0.
$$
Hence, by $\www w(0,\cdot) = 0$ in $(0,\ell)$, we obtain
$$
\www w(t,x) = \int^t_0 R(t-\tau)z(\tau,x) d\tau, \quad 0<t<T,\, x\in (0,\ell),
$$
that is,
\begin{equation}\label{(4.4)}
\www w(t,x) = \int^t_0 R(\tau)z(t-\tau,x) d\tau, \quad 0<t<T,\, x\in (0,\ell).
\end{equation}

We will prove that $z \in H^{1,2}((0,T)\times (0,\ell))$ satisfies 
\begin{equation}\label{(4.5)}
\sigma\ppp_tz(t,x) - \partial_{xx}z(t,x)+Q(x)z(t,x)=0 \quad \mbox{in}\,\, 
(0,T)\times (0,\ell)
\end{equation}
and 
\begin{equation}\label{(4.6)}
z(t,0) = \partial_x z(t,0) = 0 \quad \mbox{in $(0,T) $}.
\end{equation}

First we will verify (\ref{(4.6)}).  Indeed, $\ppp_t\www w(t,0) 
= \partial_x\ppp_t\www w(t,0) = 0$ in $(0,T)$ by (3.3), and so
$(Mz)(\cdot,0) = (M\partial_x z)(\cdot,0))= 0$ in $(0,T)$.
Consequently, the injectivity of $K_1$ directly yields (\ref{(4.6)}).

Next, we will prove (\ref{(4.5)}). 
First, using $\ppp_t\www w \in H^{1,2}((0,T)\times (0,\ell)) 
\subset C([0,T];L^2(0,\ell))$,
by (3.4) we have 
$$
\ppp_t\www w(0,x) = R(0)z(0,x), \quad x\in (0,\ell).
$$
On the other hand, substituting $t=0$ in (3.3),  we obtain
$$
\ppp_t\www w(0,x) = -\frac 1 \sigma R(0)K(x,0), \quad x\in (0,\ell).
$$
Hence $R(0)z(0,x) = -\frac 1\sigma R(0)K(x,0)$ for $x\in (0,\ell)$.  
By $R(0) \ne 0$, we reach 
\begin{equation}\label{(4.7)}
z(0,x) = -\frac 1\sigma K(x,0), \quad x\in (0,\ell).                  
\end{equation}
We proceed to the completion of the proof of (\ref{(4.5)}).  
In terms of (\ref{(4.4)}) and (\ref{(4.7)}), we have
\begin{align*}
& \sigma\ppp_t\www w(t,x) = \sigma R(t)z(0,x) 
+ \int^t_0 R(\tau)\sigma \ppp_tz(t-\tau,x) d\tau\\
=& -R(t)K(x,0) + \int^t_0 R(\tau)\sigma\ppp_tz(t-\tau,x) d\tau
\end{align*}
and
$$
(\partial_x^2-Q(x))\www w(t,x) = \int^t_0 R(\tau)
(\partial_x^2-Q(x))z(t-\tau,x) d\tau, 
\quad 0<t<T,\, x\in (0,\ell).
$$
Consequently (3.3) implies
\begin{equation}\nonumber\begin{split}
& -R(t)K(x,0) = \sigma \ppp_t\www w - \partial_x^2\www w(t,x)
+ Q\www w(t,x)\\
=& -R(t)K(x,0) + \int^t_0 R(\tau)(\sigma \ppp_tz - \partial_x^2z+Qz)
(t-\tau,x) d\tau,
\end{split}
\end{equation}
that is,
$$
\int^t_0 R(\tau)(\sigma \ppp_tz - \partial_x^2z+Qz)(t-\tau,x) d\tau = 0, \quad 
0<t<T,\, x\in (0,\ell).
$$
We set $Z(\sigma^*)(\tau):= ((\sigma\ppp_\tau z - \partial_x^2z+Qz)
(\tau,\cdot),\, \sigma^*)
_{L^2(0,\ell)}$ for any $\sigma^* \in L^2(0,\ell)$.  Then 
$$
\int^t_0 R(\tau) Z(\sigma^*)(t-\tau) d\tau = 0, \quad 0<t<T.
$$
By the Titchmarsh convolution  theorem (e.g., Titchmarsh \cite{Tit}), 
there exists $t_*(\sigma^*) \in [0,T]$ such that 
$$
R(\tau) = 0 \quad \mbox{for $0< \tau < T-t_*(\sigma^*)$}, \quad
Z(\sigma^*)(\tau) = 0 \quad \mbox{for $0 < \tau < t_*(\sigma^*)$}.
$$
If $t_*(\sigma^*) < T$, then $R(\tau) =  0$ for $0<\tau < T-t_*(\sigma^*)$ 
with $T-t_*(\sigma^*) > 0$, which implies $R(0) = 0$.  This is impossible 
because we assume (\ref{mobik}), by which we have $R(0) \ne 0$.  
Therefore $t_*(\sigma^*) = T$.
Hence, $((\sigma\ppp_{\tau}z - \partial_x^2z+Qz)(\tau,\cdot),\, \sigma^*)
_{L^2(0,\ell)} = 0$ 
for any $\tau \in (0,T)$ and any $\sigma^* \in L^2(0,\ell)$.  
Thus the verification of (\ref{(4.5)}) is complete.
\\

Now we can complete the proof of Theorem \ref{Moon1}.
Indeed, from (\ref{(4.5)}) and (\ref{(4.6)}), we can apply
the unique continuation for the parabolic equation, so that 
we reach $z=0$ in $(0,T) \times (0,\ell)$.  In particular,
$z(0,\cdot)=0$ for $0<x<\ell$.
Then equality (\ref{(4.7)}) implies $K(x,0)=0$ on $(0,\ell)$.
By the uniqueness of solution to an initial value problem 
of the hyperbolic equation in (2.1) with zero initial values
$K(x,0) = \ppp_yK(x,0) = 0$ on $\{ (x,0); 0<x<\ell\}$, we 
have $K(x,y)=0$ if $(x,y)$ is in the interior of the triangle
with vertices $(0,0)$, $(\ell,0)$ and $(\frac{\ell}{2},\, \frac{\ell}{2})$. 
In particular $K(x,x)=0$ on $(0,\frac{\ell}{2})$. This implies 
$P=Q$ in $\left(0, \,\frac{\ell}{2}\right)$.  Thus the proof of 
Theorem \ref{Moon1} is complete.
$\blacksquare$
\\
\vspace{0.2cm}
\\
{\bf Proof of Corollary \ref{Moon3}.} 
We set 
$$
R(x) = -\frac{1}{2}\int^x_0 r(y) dy, \quad 
\www R(x) = -\frac{1}{2}\int^x_0 \www r(y) dy, \quad 0<x<\ell
$$
and
$$
W(t,x) : = e^{R(x)}u(t,x), \quad \www{W}(t,x):= e^{\www R(x)}\www u.
$$
Then,  
$$
\sigma\ppp_t W = \ppp_x^2 W - \left( \frac{1}{4}r^2
- \frac{1}{2}r' + P(x)\right)W,  
\quad W(0,x) = 0, \quad 0<x<\ell, \, 0<t<T,
$$
and
$$
\sigma\ppp_t\www W = \ppp_x^2\www W
- \left( \frac{1}{4}\www r^2 - \frac{1}{2}{\www r}' + Q(x)\right)\www W,  
\quad \www W(0,x) = 0, \quad 0<x<\ell, \, 0<t<T.
$$
By (\ref{sarkofag}), we have
$$
W(t,0) = \www W(t,0)\quad \mbox{and}\quad\partial_x W(t,0) = 
\partial_x\www W(t,0)  \quad   \mbox{in}\, \, (0,T).
$$
We can readily check the conditions on $W$ and $\www W$ in Theorem \ref{Moon1},
so that the application of Theorem \ref{Moon1} completes the proof of
Corollary \ref{Moon3}.
$\blacksquare$

\section{ Proof of Theorem \ref{AX}}\label{four}

Let $\ep \in (0,\ell]$ be the maximal number satisfying 
\begin{equation}\label{xui1} 
P(x)=Q(x)\quad \mbox{in}\quad (0,\ep).
\end{equation}
By the assumption of the theorem, such $\ep>0$ exists.
We can assume $\ep < \ell$. Otherwise, the theorem is already proved. 
We choose $x_0>0$ sufficiently small such that $0<x_0<\ep$.

We introduce a transformation of the variables: $x \mapsto \eta$ by 
$\eta = x-\ep + x_0$ for $\ep - x_0 < x < \ell$.
We note that $\ep - x_0 < x < \ell$ if and only if
$0< \eta < \ell + x_0 - \ep$, and in particular,
$\ep - x_0 < x < \ep$ if and only $0<\eta<x_0$.
We set $P(\eta;x_0) := P(x)$,  $Q(\eta;x_0) := Q(x)$, 
$u^k(t,\eta;x_0) := u^k(t,x)$ and $\www u^k(t,\eta;x_0) := \www u^k(t,x)$ for 
$1\le k \le N$.

Then, (\ref{xui1}) implies
\begin{equation}\label{mkP}
P(\eta;x_0) = Q(\eta;x_0) \quad \mbox{for $0<\eta<x_0$}.   
\end{equation}

We set $\OOO_0:= \{ (\eta,y);\, 0<y<x<\ell-\ep+x_0\}$.
In terms of Lemma 1, there exists a unique solution
$K(\eta,y;x_0) \in C^2(\ooo{\OOO_0})$ to 
\begin{equation}\label{(7.1)}
\left\{ \begin{array}{rl}
& \ppp_\eta^2 K(\eta,y;x_0) - \ppp_y^2K(\eta,y:x_0) 
= Q(\eta;x_0)\KKK - \KKK P(y;x_0) \quad
\mbox{in $\OOO_0$},\\ 
& \ppp_y K(\eta,0;x_0) = 0, \quad 0<\eta<\ell-\ep+x_0, \\
& 2\frac{d}{d\eta}K(\eta,\eta;x_0) = Q(\eta;x_0) - P(\eta;x_0),
\quad 0<\eta<\ell-\ep+x_0, \quad  K(0,0; x_0) = 0.
\end{array}\right.   
\end{equation}  

Setting 
$$
V^k(t,\eta):= u^k(t,\eta;x_0) + \int^{\eta}_0 K(\eta,y;x_0)u^k(t,y;x_0)dy,
\quad 1\le k \le N,
$$
by Lemma \ref{main} we obtain
\begin{equation}\label{XMmk1}
\left\{ \begin{array}{rl}
&\sigma\ppp_t V^k - \ppp_\eta^2 V^k + Q(\eta;x_0)V^k 
= -K(\eta,0;x_0)\partial_\eta u^k(t,0;x_0)\quad
\mbox{in $(0,T)\times (0,\ell-\ep + x_0)$},          \\
&\partial_\eta V^k(t,0)=\partial_\eta {u}^k(t,0;x_0)
   \quad \mbox{and}\quad  V^k(t,0;x_0) = u^k(t,0;x_0) 
\quad \mbox{in $(0,T)$}.
\end{array}\right.
\end{equation}

By (4.2), we see
$$
K(\eta,\eta;x_0) = 0 \quad \mbox{for $0<\eta<x_0$}.
$$
Therefore, the uniqueness of solution to problem (4.3) yields
\begin{equation}\label{4.5}
K(\eta,y;x_0) = 0, \quad  0 < y < \eta < x_0. 
\end{equation}
Hence, substituting $y=0$, we reach 
$$
K(\eta,0;x_0) = 0 \quad \mbox{for $0<\eta<x_0$}.
$$
With (4.3), we obtain
\begin{equation}\label{(P3.1)}
\left\{ \begin{array}{rl}
& \ppp_\eta^2K(\eta,y;x_0) - \ppp_y^2K(\eta,y;x_0) 
= Q(\eta;x_0)K(\eta,y;x_0) - K(\eta,y;x_0)P(y;x_0),\\
& \quad \qquad \quad \qquad\qquad \qquad \qquad \mbox{for $0<y<\eta<x_0$}, \\
& K(\eta,0;x_0) = \ppp_y K(\eta,0;x_0) = 0, \quad 0<\eta<x_0.
\end{array}\right.
\end{equation}
Moreover, we make the even extension with respect to $y=0$ to
$\OOO_{\pm}:= \{ (\eta,y);\, -\eta < y < \eta, \, 
0<\eta < x_0+\ell-\ep\}$:
$$
K(\eta,y:x_0) = K(\eta,-y;x_0) \quad 
\mbox{for $-\eta < y < 0, \, 0<\eta<x_0 + \ell-\ep$}.
$$
In view of $\ppp_yK(\eta,0;x_0) = 0$ for $0< \eta < x_0+\ell-\ep$, we can 
readily verify that the extended $K(\eta,y;x_0)$ is in $C^2(\ooo{\OOO_{\pm}})$,
and satisfies the hyperbolic equation in (4.3) in $\OOO_{\pm}$.
Therefore, setting $\delta:= \min\{ x_0, \, \ell-\ep\}$, we have
\begin{equation}\label{4.7}
\left\{ \begin{array}{rl}
& \ppp_\eta^2K(\eta,y;x_0) - \ppp_y^2K(\eta,y;x_0) 
= Q(\eta;x_0)K(\eta,y;x_0) - K(\eta,y;x_0)P(y;x_0), \\
& K(x_0,y:x_0) = \ppp_{\eta}K(x_0,y;x_0) = 0 \quad
\mbox{for $-x_0 < y < x_0$.}
\end{array}\right.
\end{equation}
The uniqueness of solution to the hyperbolic equation with 
initial condition (4.7) implies that 
$K(\eta,y;x_0) = 0$ if $x_0 < \eta < x_0 + \delta$ and 
$\eta - 2x_0 < y < -\eta + 2x_0$, from which we see
$K(\eta,0;x_0) = 0$ for $x_0 < \eta < x_0+\delta$ with $\delta > 0$.
Then 
$$
K(\eta,0;x_0)\ppp_xu^k(t,0) = 0 \quad
\mbox{for $k\in \{1,2,..., N\}$, $0<\eta<x_0+\delta$,
$0<t<T$}
$$
and $P(\eta;x_0) = Q(\eta;x_0)$ for $0<\eta<x_0$.
In view of (1.3) - (1.4) - (1.5), we can apply Lemma 3 and we can find 
a constant $\ep_1 > 0$ such that $P(\eta;x_0) = Q(\eta;x_0)$ for $0<\eta<
x_0+\ep_1$, that is, $P(x) = Q(x)$ for $0<x<\ep_1+ \ep$.
This is a contradiction against the choice of $\ep$ in (4.1).
Thus the proof of Theorem 6 is complete.
$\blacksquare$

\section{Concluding remarks}\label{five}
\begin{itemize}
\item
We can describe our main achievements:
\\
(i) Uniqueness without full boundary conditions in spatially one-dimensional
cases.
\\
(ii) Wide applicability of our method which is based on 
the transformation operator but does not depend on any results on 
inverse spectral problems.
\item
In this article, we mainly consider equations of the form 
$$
\sigma\ppp_tu - \ppp_x^2u + P(x)u = 0,      \eqno{(5.1)}
$$
but by means of the classical Liouville
transform, we can reduce an inverse problem of determining 
$p(x)$ in 
$$
\sigma\ppp_tu(t,x) = \ppp_x(p(x)\ppp_x u(t,x))
$$
to the inverse problem for (5.1).  We omit the details.
\item
Inverse parabolic problems with initial values are
difficult and the uniqueness is not known in general.
In order to solve it, one way is to change the inverse 
parabolic problem to an inverse hyperbolic problem.

One traditional way is that 
by considering an integral transform in time of solutions to 
parabolic equations, we can reduce the inverse parabolic
problem to an inverse hyperbolic problem where the 
uniqueness is classical with initial values.
Such an integral transform is similar to the Laplace transform and 
called Reznitskaya's transform in the context of the inverse problem
(e.g., Romanov \cite{Ro}),   
but we need data of the solutions over the time interval
$(0, \infty)$.  In the case where the solution data are time analytic,
we can reduce data over $(0,\infty)$ to data over a finite time interval and 
the analyticity is valid if the boundary values are time analytic and
the coefficients of the parabolic equation are independent of the time 
or analytic in time.  As for this approach, see \cite{IY230}, \cite{Kli}.
However, in our case, since we do not assume the boundary 
values on the whole boundary $x=0$ and $x=\ell$, we can not expect 
the time analyticity of $u$ and $\www{u}$.
\end{itemize}

{\bf Acknowledgments.}
The work was supported by Grant-in-Aid for Scientific Research (A) 20H00117 
and Grant-in-Aid for Challenging Research (Pioneering) 21K18142 of 
Japan Society for the Promotion of Science.

\end{document}